\documentclass[12pt]{article}
\usepackage{amsmath}
\usepackage{amsfonts}
\DeclareRobustCommand{\qbinom}{\genfrac[]{0pt}{}}

\newtheorem{theorem}{Theorem}
\newtheorem{lemma}{Lemma}

\begin{document}
\title{Little $q$-Legendre polynomials and irrationality of certain
Lambert series}
\author{Walter Van Assche
\thanks{WVA is a Research Director of the Belgian Fund for Scientific Research 
(FWO).}
\\
\normalsize{Katholieke Universiteit Leuven and Georgia Institute of Technology}}
\date{\today}
\maketitle

\begin{abstract}
Certain $q$-analogs $h_p(1)$ of the harmonic series, with $p=1/q$ an integer 
greater than one, were shown to be irrational by
Erd\H{o}s \cite{erdos}. In 1991--1992 Peter Borwein \cite{bor1} \cite{bor2}
 used Pad\'e approximation
and complex analysis to prove the irrationality of these $q$-harmonic series
and of $q$-analogs $\ln_p(2)$ of the natural logarithm of 2. Recently
Amdeberhan and Zeilberger \cite{WZ} used the qEKHAD symbolic package to find
$q$-WZ pairs that provide a proof of irrationality similar to Ap\'ery's
proof of irrationality of $\zeta(2)$ and $\zeta(3)$. They also obtain
an upper bound for the measure of irrationality, but  better upper
bounds were earlier given by Bundschuh and V\"a\"an\"anen \cite{BV} and recently
also by Matala-aho and V\"a\"an\"anen \cite{MV} (for $\ln_p(2)$). 
In this paper we show how one can obtain
rational approximants for $h_p(1)$ and $\ln_p(2)$ (and many other
similar quantities) by  Pad\'e approximation using
 little $q$-Legendre polynomials
and we show that properties of these orthogonal polynomials indeed
prove the irrationality, with an upper bound of the measure
of irrationality  which is as sharp as the upper bound given by
Bundschuh and V\"a\"an\"anen for $h_p(1)$ and a
better upper bound as the one given by Matala-aho and V\"a\"an\"anen
for $\ln_p(2)$.
\end{abstract}

\section{Introduction}

Most important special functions, in particular hypergeometric functions,
have $q$-ex\-ten\-sions, usually obtained by replacing Pochhammer symbols
$(a)_n = a(a+1)\cdots(a+n-1)$ by their $q$-analog
$(a;q)_n = (1-a)(1-aq)(1-aq^2)\cdots(1-aq^{n-1})$. Since
\[   \lim_{q \to 1} \frac{(q^a;q)_n}{(1-q)^n} = (a)_n,  \]
one usually retrieves the original special function from its
$q$-extension by letting $q \to 1$. A good source for
$q$-extensions of hypergeometric series (basic hypergeometric
series) is the book \cite{GR} by Gasper and Rahman. Our interest
in this paper is the following $q$-extension of the harmonic series
\begin{equation}  \label{eq:hp1}
  h_p(1) = \sum_{k=1}^\infty \frac{1}{p^k-1} 
        = \sum_{k=1}^\infty \frac{q^k}{1-q^k}, 
        \qquad 0 < q = \frac1p < 1, 
\end{equation}
 and of the natural logarithm of 2
\begin{equation}  \label{eq:logp2}
    \ln_p(2) = \sum_{k=1}^\infty  \frac{(-1)^k}{p^k-1} 
             = \sum_{k=1}^\infty \frac{(-q)^k}{1-q^k},
               \qquad 0 < q = \frac1p < 1.
\end{equation}               
In 1948, Paul Erd\H{o}s \cite{erdos} proved that $h_2(1)$ is irrational.
Peter Borwein \cite{bor1} \cite{bor2} showed that $h_p(1)$ 
(and other similar numbers) are
irrational for every integer $p >1$ and also proved the
irrationality of $\ln_p(2)$ for every integer $p > 1$. If we denote 
$E_p(z) = (qz;q)_\infty$, with $q=1/p$, then
B\'ezivin \cite{Be} had earlier shown that $E_p(\alpha), E_p'(\alpha), \ldots, E_p^{(k)}(\alpha)$ are linearly independent over $\mathbb{Q}$
for every $k \in \mathbb{N}$ and $\alpha \in \mathbb{Q}$ with $\alpha p^j
\neq -1$ for every integer $j$. The case $k=1$ corresponds to irrationality
of
\[   \sum_{j=1}^\infty \frac{1}{\alpha+p^j}. \]
Borwein used
Pad\'e approximation techniques and complex analysis to obtain good rational
approximants to $h_p(1)$ and $\ln_p(2)$. Indeed, one can use the following lemma
to prove irrationality \cite[Lemma 5.1]{wva}:
\begin{lemma} \label{lem1}
Let $x$ be a
real number, and suppose there exist integers $a_n, b_n$ 
$(n \in \mathbb{N})$ such that
\begin{enumerate}
 \item $x\neq a_n/b_n$ for every $n \in \mathbb{N}$,
 \item $\lim_{n \to \infty} ( b_n x-a_n) = 0$,
\end{enumerate}
then $x$ will be irrational.
\end{lemma}
This lemma expresses the fact that the order of approximation of a
rational number by rational numbers is one and not higher \cite[Theorem 186]{HW}.
Furthermore,
if $|x-a_n/b_n| = \mathcal{O}(1/b_n^{1+s})$ with $0 < s < 1$
and $b_n < b_{n+1} < b_n^{1+o(1)}$, then the
measure of irrationality $r(x)$ (Liouville-Roth number, order of approximation)
\[ r(x) = \inf \{ r \in \mathbb{R}: |x-a/b| < 1/b^r \textrm{\ has\ at\ most\
 finitely\ many\ integer\ solutions\ } (a,b) \} \]
 satisfies
$2 \leq r(x) \leq 1+1/s$ (see, e.g., \cite[exercise 3 on p.~376]{BB}
for the upper bound; the lower bound follows since every irrational number is approximable
to order 2 \cite[Theorem 187]{HW}).   Recently Amdeberhan and Zeilberger
\cite{WZ} found $q$-WZ pairs to obtain rational approximants
for $h_p(1)$ and $\ln_p(2)$, improving the upper bound
for the measure of irrationality to $4.8 = 24/5$. However,
four years earlier Bundschuh and V\"a\"an\"anen \cite[p.~178]{BV}
had established better upper bounds: $r(h_p(1)) \leq 1+(\pi^2+2)/(\pi^2-2)
=2.50828\ldots$ and $r(\ln_p(2)) \leq 1+(2\pi^2+3)/(\pi^2-3)=4.310119\ldots$.
The upper bound for $r(\ln_p(2))$ was improved by Matala-aho and
V\"a\"an\"anen \cite{MV} to $3.9461$.

In this paper
we show that one can find rational approximants which are related to
little $q$-Legendre polynomials and hence return to Pad\'e approximation.
We can then use some results for little $q$-Legendre polynomials
to prove the irrationality once more, and with the aid of some
elementary number theory we obtain the same bound for the measure
of irrationality as the one obtained by Bundschuh and V\"a\"an\"anen
for $h_p(1)$ and a better upper bound as the one given by Matala-aho 
and V\"a\"an\"anen for $\ln_p(2)$.
The connection with little $q$-Legendre polynomials opens the
way for proving the irrationality of $q$-extensions of
$\zeta(2)$ and $\zeta(3)$ in Ap\'ery's spirit \cite{apery},
 using multiple orthogonal $q$-polynomials \cite{wva}.
  
\section{Little $q$-Legendre polynomials}

The little $q$-Legendre polynomials are defined by
\begin{eqnarray}   \label{eq:Pnqphi}
  P_n(x|q) & = & {}_2\phi_1 \left(  \left. 
     \begin{array}{c} q^{-n} , q^{n+1} \\ q \end{array}
      \right| q ; qx \right)   \nonumber \\
      & = & \sum_{k=0}^n \frac{ (q^{-n};q)_k (q^{n+1};q)_k}{(q;q)_k}
        \frac{q^k x^k}{(q;q)_k}, \qquad  0 < q < 1,
\end{eqnarray}    
and they are orthogonal polynomials on the exponential lattice
$\{ q^k, k=0,1,2,\ldots\}$:
\begin{equation}   \label{eq:Pnorth}
  \sum_{k=0}^\infty  q^k P_m(q^k|q) P_n(q^k|q) = 
  \frac{q^n}{1-q^{2n+1}} \delta_{n,m}. 
\end{equation}
If we use the $q$-binomial coefficients
\begin{equation}  \label{eq:qbinom}
  \qbinom{n}{k}_q = \frac{(q;q)_n}{(q;q)_k(q;q)_{n-k}} \ ,
\end{equation}
and the formulas
\begin{equation}  \label{eq:qminplus}
  (q^{-n};q)_k = (-1)^k q^{-nk+k(k-1)/2} \frac{(q;q)_n}{(q;q)_{n-k}}, 
  \quad 
  (q^{n+1};q)_k = \frac{(q;q)_{n+k}}{(q;q)_n}, 
\end{equation}
then (\ref{eq:Pnqphi}) reduces to
\begin{equation}  \label{eq:Pnqbinom}
  P_n(x|q) = \sum_{k=0}^n \qbinom{n}{k}_q \qbinom{n+k}{k}_q 
   q^{-nk+k(k+1)/2} (-x)^k,
\end{equation}     
and since 
\[  \lim_{q \uparrow 1}  \qbinom{n}{k}_q = \binom{n}{k}, \]
we see that we indeed find the Legendre polynomials on
$[0,1]$ by letting $q$ tend to 1
\[    \lim_{q \uparrow 1} P_n(x|q) = 
   \sum_{k=0}^n \binom{n}{k} \binom{n+k}{k}  (-x)^k = P_n(x). \]
Let $p = 1/q$ so that $p > 1$ whenever $0 < q < 1$, then
\begin{equation} \label{eq:qtop}
   (q;q)_k = (-1)^k p^{-k(k+1)/2} (p;p)_k, \quad
   \qbinom{n}{k}_q = p^{-k(n-k)} \qbinom{n}{k}_p,
\end{equation}
and we can rewrite the little $q$-Legendre polynomial as
\begin{equation}  \label{eq:Pnpbinom}
  P_n(x|q) = \sum_{k=0}^n \qbinom{n}{k}_p \qbinom{n+k}{k}_p
    p^{-kn+k(k-1)/2} (-x)^k.
\end{equation}
There is a Rodrigues formula for these little $q$-Legendre
polynomials in terms of the $q$-difference operator $D_q$
for which
\begin{equation}  \label{eq:Dq}
   D_qf(z) = \begin{cases}
              {\displaystyle \frac{f(z)-f(qz)}{(1-q)z}} & \textrm{if } z\neq 0, \\
               f'(0) & \textrm{if } z = 0,
              \end{cases}
 \end{equation}                       
namely
\begin{equation}  \label{eq:qRod}
  P_n(x|q) = \frac{q^{n(n-1)/2} (1-q)^n}{(q;q)_n} D_p^n [ (qx;q)_n x^n ],
\end{equation}
which will be useful later. We refer to \cite{koek} for more information
and references for little $q$-Legendre polynomials.

Equation (\ref{eq:Pnpbinom}) expresses the little $q$-Legendre polynomials
in the basis $\{1,x,x^2,\ldots,x^n\}$ of monomials. Sometimes it is more convenient
to use another basis of polynomials, and for orthogonal polynomials  on 
$\{q^k,\ k=0,1,2,\ldots\}$ a convenient set of basis functions
is $\{(qx;q)_k,\ k=0,1,2,\ldots,n\}$. We will need to use some
$q$-series for this purpose. Recall the $q$-analog of Newton's binomium formula
\begin{equation}  \label{eq:qbinomium}
  (x;q)_n = \sum_{k=0}^n \qbinom{n}{k}_q q^{k(k-1)/2} (-x)^k ,
\end{equation}
and its dual
\begin{equation}  \label{eq:qbinomium2}
  x^n = \sum_{k=0}^n \qbinom{n}{k}_q (-1)^k q^{-nk+k(k+1)/2} (x;q)_k.   
\end{equation}
A more general result is the $q$-binomial series
\begin{equation}  \label{eq:qbinomial}
  \sum_{n=0}^\infty \frac{(a;q)_n}{(q;q)_n} x^n = \frac{(ax;q)_\infty}{(x;q)_\infty}, \qquad
  |q| < 1, |x| < 1.
\end{equation}
Use (\ref{eq:qbinomium2}) with argument $q^{n+1}x$ in  the Rodrigues
formula (\ref{eq:qRod}), then we have
\[ P_n(x|q) = \frac{p^{n(n+1)}(1-p)^n}{(p;p)_n} 
  \sum_{k=0}^n \qbinom{n}{k}_q (-1)^k q^{-nk+k(k+1)/2}
   D_p^n (qx;q)_{n+k} . \]
One finds easily that
\begin{equation}  \label{eq:Dppoch}
  D_p^k (qx;q)_n = \frac{(q;q)_n}{(q;q)_{n-k}(1-p)^k} (qx;q)_{n-k} ,
\end{equation}      
so that, using (\ref{eq:qtop}) we find the required expansion 
\begin{equation}  \label{eq:Pnpoch}
 P_n(x|q) = (-1)^n  \sum_{k=0}^n \qbinom{n}{k}_p \qbinom{n+k}{k}_p (-1)^k p^{(n-k)(n-k+1)/2} (qx;q)_k.
\end{equation}

\section{The $q$-harmonic series}

Orthogonal polynomials naturally arise in Pad\'e approximation of a
Stieltjes function
\begin{equation}  \label{eq:markov} 
  f(z) = \int \frac{d\mu(x)}{z-x}, \qquad z \notin \textrm{supp}(\mu).
\end{equation}
Suppose $\mu$ is a positive measure on the real line with infinite
support and for which all the moments exist. If $P_n(x)$ $(n = 0,1,2,\ldots)$
are orthogonal polynomials for $\mu$, i.e., $P_n$ is of degree $n$ and
\[  \int P_n(x)P_m(x)\, d\mu(x) = 0, \qquad m\neq n,  \]
and if $Q_n$ are the polynomials of degree $n-1$ given by
\begin{equation}  \label{eq:Qn}
  Q_n(z) = \int \frac{P_n(z)-P_n(x)}{z-x} \, d\mu(x), 
\end{equation}
then it is easy to see that
\begin{equation}  \label{eq:pade}
  P_n(z) f(z) - Q_n(z) = \int \frac{P_n(x)}{z-x} \, d\mu(x), 
  \qquad z \notin \textrm{supp}(\mu).   
\end{equation}        
If we expand $1/(z-x)$ around $z=\infty$ as
\[   \frac{1}{z-x} = \sum_{k=0}^{n-1} \frac{x^k}{z^{k+1}}
    + \frac{x^n}{z^n} \frac{1}{z-x} , \]
then by orthogonality we have
\[   P_n(z) f(z) - Q_n(z) = \frac{1}{z^n} \int \frac{P_n(x) x^n}{z-x}\, d\mu(x)
    = \mathcal{O}(1/z^{n+1}). \]
This is precisely the (linearized) form of the interpolation conditions near
$z=\infty$ for Pad\'e approximation so that $Q_n(z)/P_n(z)$ is the 
$[\frac{n-1}{n}]$ Pad\'e
approximant for $f(z)$ near $z=\infty$. 

For little $q$-Legendre polynomials the measure $\mu$ is supported on
$\{q^k, k=0,1,2,\ldots\}$, which is a bounded set in $[0,1]$ with one accumulation
point at $0$. The measure is given by
\[  \int g(x) \, d\mu(x) = \sum_{k=0}^\infty g(q^k) q^k. \]
The Stieltjes function for this measure is
\begin{equation}  \label{eq:markov2}
  f(z) = \sum_{k=0}^\infty \frac{q^k}{z-q^k} = \sum_{k=0}^\infty \frac{1}{zp^k-1}. 
\end{equation}
We will need this function at $p^n$, where it gives
\begin{equation}  \label{eq:markovpn}
   f(p^n) = \sum_{k=0}^\infty \frac{1}{p^{n+k}-1} = 
   h_p(1) - \sum_{k=1}^{n-1} \frac{1}{p^k-1}.
\end{equation}     
Hence if $p>1$ is an integer, then $f(p^n)$ gives $h_p(1)$ up to 
$\sum_{k=1}^{n-1} 1/(p^k-1)$, which is a rational number. Now use (\ref{eq:pade}) for
little $q$-Legendre polynomials at $z=p^n$ to find
\begin{equation}  \label{eq:padepn}
  P_n(p^n|q) \left( h_p(1) - \sum_{k=1}^{n-1} \frac{1}{p^k-1} \right)
   - Q_n(p^n|q) = \sum_{k=0}^\infty \frac{P_n(q^k|q)}{p^n-q^k} q^k.
\end{equation}
Observe that (\ref{eq:Pnpbinom}) gives
\begin{equation}  \label{eq:Pnpn}
  P_n(p^n|q) = \sum_{k=0}^n \qbinom{n}{k}_p \qbinom{n+k}{n}_p (-1)^k p^{k(k-1)/2},
\end{equation}
which is nearly the $b_n$ found in \cite[p.~277]{WZ} (their $b_n$ corresponds
to $P_n(p^{n+1}|q)$). Observe also that Borwein's  
construction \cite[Lemma 2]{bor2} uses $P_{n-1}(cp^{n+1}|q)$.
The $p$-binomial numbers $\qbinom{n}{k}_p$ are polynomials in $p$ with
integer coefficients, which follows easily from the $q$-version
of Pascal's triangle identities
\[  \qbinom{n}{k}_p = \qbinom{n-1}{k-1}_p + p^k \qbinom{n-1}{k}_p =
    \qbinom{n-1}{k}_p + p^{n-k} \qbinom{n-1}{k-1}_p, \]
hence if $p >1$ is an integer then $\qbinom{n}{k}_p$ and
$\qbinom{n+k}{k}_p$ are integers. This means that (\ref{eq:Pnpn}) implies 
$P_n(p^n|q)$ to be an integer. Furthermore, since $p^n > 1$ and all the
zeros of $P_n(x|q)$ are in $[0,1]$, we also may conclude
that $(-1)^n P_n(p^n|q)$ is positive for all $n$.

The second important quantity in (\ref{eq:padepn}) is $Q_n(p^n|q)$. This
associated little $q$-Legendre polynomial can be computed explicitly using
(\ref{eq:Qn}) and is given by
\[  Q_n(x|q) = \sum_{j=0}^\infty \frac{P_n(x|q)-P_n(q^j|q)}{x-q^j} q^j. \]
Use (\ref{eq:Pnpoch}) to write this as
\[  Q_n(x|q) = (-1)^n \sum_{k=0}^{n} \qbinom{n}{k}_p \qbinom{n+k}{k}_p (-1)^k
   p^{(n-k)(n-k+1)/2}   \sum_{j=0}^\infty \frac{(qx;q)_k-(q^{j+1};q)_k}{x-q^j} q^j. \]
Now use
\[   \frac{(qx;q)_k - (qy;q)_k}{x-y} = - \sum_{\ell=1}^k q^\ell (qy;q)_{\ell-1}
(q^{\ell+1}x;q)_{k-\ell} , \]
which one can easily prove by induction,
then this gives
\begin{multline*}
 Q_n(x|q) =  (-1)^{n+1} \sum_{k=0}^n \qbinom{n}{k}_p \qbinom{n+k}{k}_p (-1)^k
                  p^{(n-k)(n-k+1)/2}  \\
           \times       \sum_{\ell=1}^{k} q^\ell (q^{\ell+1}x;q)_{k-\ell}
                   \sum_{j=0}^\infty q^j (q^{j+1};q)_{\ell-1}. 
\end{multline*}
Using the $q$-binomial series (\ref{eq:qbinomial}) we
can calculate the modified moments
\begin{eqnarray*} 
  \sum_{j=0}^\infty q^j (q^{j+1};q)_{\ell-1} & = & (q;q)_{\ell-1}
   \sum_{j=0}^\infty q^j \frac{(q^\ell;q)_j}{(q;q)_j} \\
  & = & \frac{(q;q)_{\ell-1} (q^{\ell+1};q)_\infty}{(q;q)_\infty} \\
  & = & \frac{1}{1-q^\ell}, 
\end{eqnarray*}
so that
\begin{equation} \label{eq:Qnx}
  Q_n(x|q) = (-1)^{n+1} \sum_{k=0}^n \qbinom{n}{k}_p \qbinom{n+k}{k}_p (-1)^k
                  p^{(n-k)(n-k+1)/2} \sum_{\ell=1}^{k} \frac{(q^{\ell+1}x;q)_{k-\ell}}{p^\ell-1} . 
\end{equation}                            
Evaluating at $x=p^n$, and using 
\[  (q^{\ell+1}p^n;q)_{k-\ell} = (p^{n-k};p)_{k-\ell}, \]
then gives
\begin{equation}  \label{eq:Qnpn}
  Q_n(p^n|q) = (-1)^{n+1} \sum_{k=0}^n \qbinom{n}{k}_p \qbinom{n+k}{k}_p (-1)^k p^{(n-k)(n-k+1)/2}
                   \sum_{\ell=1}^{k} \frac{(p^{n-k};p)_{k-\ell}}{p^\ell-1}.
\end{equation}
 All the terms
in the sum for $Q_n(p^n|q)$ are now integers, except
for the $p^{\ell}-1$ in the denominators. In order to obtain an integer
we therefore need to multiply everything by a multiple of
all $p^{\ell}-1$ for $\ell=1,2,\ldots,n$. 
We will choose
\begin{equation}  \label{eq:lcm}
  d_n(p) = \prod_{k=1}^n \Phi_k(p),
\end{equation}  
where
\begin{equation}
    \Phi_n(x) = \prod_{\genfrac{}{}{0pt}1{k=1}{\gcd(k,n)=1}}^n 
   (x-\omega_n^k), 
\qquad \omega_n = e^{2\pi i/n},
\end{equation}
are the cyclotomic polynomials \cite[\S 4.8]{St}. Each cyclotomic polynomial
is monic, has integer coefficients, and the degree of
$\Phi_n$ is $\phi(n)$ (Euler's totient function). It is known that
\begin{equation}   \label{eq:xnPhi}
   x^n-1 = \prod_{d | n} \Phi_d(x), 
\end{equation}
and that every cyclotomic polynomial is irreducible over $\mathbb{Q}[x]$.
Hence $d_n(p)$ is a multiple of all $p^\ell-1$ for $\ell=1,2,\ldots,n$.
The growth of this sequence is given by the following lemma, which
was essentially given by A. O. Gel'fond, who obtained
the upper bound in \cite[Equation (7)]{Ge}. We give a proof to
make this paper self-contained.

\begin{lemma}  \label{lm:dn}
Suppose $p$ is an integer greater than one and let $d_n(p)$
be given by (\ref{eq:lcm}). Then
\begin{equation}  \label{eq:dn}
    \lim_{n \to \infty} d_n(p)^{1/n^2} = p^{3/\pi^2}.
\end{equation}
\end{lemma}
\textbf{Proof:}
The degree of $d_n(p)$ as a monic polynomial in $p$ is $\sum_{k=1}^n \varphi(k)$
and an old result of Mertens (1874) shows that for $n \to \infty$ this 
grows like \cite[Theorem 330]{HW}
\begin{equation}   \label{eq:mertens}
  \sum_{k=1}^n \varphi(k) = \frac{3}{\pi^2} n^2 + \mathcal{O}(n\log n). 
\end{equation}
We will adapt the classical proof of (\ref{eq:mertens}) to prove our lemma.
For this we use M\"obius inversion of (\ref{eq:xnPhi}) to find the representation
\[   \Phi_n(x) = \prod_{d | n} (x^d-1)^{\mu(n/d)}, \]
where $\mu$ is the M\"obius function. Taking logarithms in (\ref{eq:lcm})
gives
\[    \log d_n(p) = \sum_{k=1}^n \sum_{d | k} \mu(k/d) \log (p^d-1). \]
Changing the order of summation gives
\[      \log d_n(p) = \sum_{\ell=1}^n \mu(\ell) \sum_{d=1}^{\lfloor n/\ell \rfloor} \log (p^d-1). \]
Now
\[   \sum_{d=1}^m \log (p^d-1) = \log  p^{m(m+1)/2} (q;q)_m = \frac{m(m+1)}2 \log p + \log (q;q)_m, \]
so that
\begin{eqnarray*}
  \log d_n(p) & = & \frac{n^2}{2} \log p \sum_{\ell=1}^n  \frac{\mu(\ell)}{\ell^2}
    + \mathcal{O}\left(n \sum_{\ell=1}^n \frac{1}{\ell} \right) \\
    & = & \frac{n^2 \log p}{2 \zeta(2)} + \mathcal{O}(n\log n). 
\end{eqnarray*} 
The lemma now follows by using $\zeta(2) = \pi^2/6$.
\ \framebox{}
\bigskip

So far we have found that the numbers
\begin{eqnarray}
   b_n & = &  d_n(p) P_n(p^n|q), \label{eq:bn} \\
     a_n & = &  d_n(p) Q_n(p^n|q) +
      b_n \sum_{k=1}^{n-1} \frac{1}{p^k-1} .       \label{eq:an}
\end{eqnarray}
are integers and         
\begin{equation}
   b_n h_p(1) - a_n =  d_n(p) \sum_{k=0}^\infty
      \frac{P_n(q^k|q)}{p^n-q^k} q^k.
\end{equation}
We now want to show that $b_n h_p(1)-a_n \neq 0$ for all $n$ and
\[  \lim_{n \to \infty} ( b_n h_p(1) - a_n ) = 0, \]
so that Lemma \ref{lem1} implies the irrationality of $h_p(1)$.
First observe that
\[  \sum_{k=0}^\infty
      \frac{P_n(q^k|q)}{p^n-q^k} q^k = \frac{1}{P_n(p^n|q)}
      \sum_{k=0}^\infty
      \frac{P_n(q^k|q)P_n(p^n|q)}{p^n-q^k} q^k. \]
 If we add and subtract $P_n(q^k|q)$ in the sum, then we have
\[  \sum_{k=0}^\infty
      \frac{P_n(q^k|q)}{p^n-q^k} q^k =
       \frac{1}{P_n(p^n|q)}
      \sum_{k=0}^\infty P_n(q^k|q) \frac{P_n(p^n|q)-P_n(q^k|q)}{p^n-q^k} q^k
+ \frac{1}{P_n(p^n|q)}
      \sum_{k=0}^\infty
      \frac{P_n^2(q^k|q)}{p^n-q^k} q^k. \]
The first sum on the right hand side vanishes because 
of the orthogonality, so that      
\begin{equation}  \label{eq:error}
   b_n h_p(1) - a_n = \frac{d_n(p)}{P_n(p^n|q)} \sum_{k=0}^\infty
      \frac{P_n^2(q^k|q)}{p^n-q^k} q^k. 
\end{equation}
All the terms in the sum are now positive, and $(-1)^n P_n(p^n|q)$
is positive for all $n$, hence we may conclude that
\[   (-1)^n (b_n h_p(1)-a_n) > 0, \qquad n=1,2,\ldots \]
Next we show that this quantity converges to zero. Clearly
$p^n-1 \leq p^n-q^k \leq p^n$ so that
\[     \frac{1}{p^n} \sum_{k=0}^\infty
      P_n^2(q^k|q) q^k \leq
\sum_{k=0}^\infty \frac{P_n^2(q^k|q)}{p^n-q^k} q^k \leq
      \frac{1}{p^n-1} \sum_{k=0}^\infty
      P_n^2(q^k|q) q^k . \]
Now we can use the norm of the little $q$-Legendre polynomial (\ref{eq:Pnorth}) to find
\begin{equation}  \label{eq:bounds}
 \frac{p}{p^{2n+1}-1} \leq
\sum_{k=0}^\infty \frac{P_n^2(q^k|q)}{p^n-q^k} q^k \leq
  \frac{p^{n+1}}{(p^n-1)(p^{2n+1}-1)}. 
\end{equation}
 What remains is to find the
 asymptotic behavior of $P_n(p^n|q)$ as $n \to \infty$. For this we can
use a very general theorem for sequences of polynomials with uniformly
bounded zeros.

\begin{lemma}  \label{lem:Pn}
Suppose $P_n$ $(n \in \mathbb{N})$ is a sequence of monic polynomials of degree $n$ and that the zeros $x_{j,n}$ $(1 \leq j \leq n)$ of $P_n$ are
such that $|x_{j,n}| \leq M$, with $M$ independent of $n$. Then
we have for $|x|>1$ and every
$c \in \mathbb{C}$
\begin{equation}       
  \lim_{n\to\infty} |P_n(cx^n)|^{1/n^2} =  |x|.
\end{equation}
\end{lemma}
\textbf{Proof:}
Factoring the polynomial $P_n$ gives
\[  |P_n(x)| = \prod_{j=1}^n |x-x_{j,n}|. \]
We have the obvious bounds
\[   |x| - M    \leq |x-x_{j,n}| \leq |x| + M, \]
hence when $|cx^n| > M$ 
\[  ( |cx^n|-M )^n \leq |P_n(cx^n)| \leq (|cx^n| + M)^n. \]
For $n$ large enough this easily gives
\[   |x| \left(|c| - \frac{M}{|x|^n} \right)^{1/n} \leq |P_n(cx^n)|^{1/n^2}
  \leq |x| \left( |c| + \frac{M}{|x|^n} \right)^{1/n}, \]
which gives the desired result.
\ \framebox{}
\bigskip

Observe that we may allow $M$ to grow with $n$  subexponentially. For
little $q$-Legendre polynomials the zeros are all in $[0,1]$ so that
we can use the Lemma with $M=1$. The
leading coefficient $\kappa_n$ of $P_n(x|q)$ is, by (\ref{eq:Pnpbinom}), equal to
$\kappa_n=(-1)^n \qbinom{2n}{n}_p\, p^{-n(n+1)/2}$, giving
\[  \lim_{n \to \infty} |\kappa_n|^{1/n^2} = \sqrt{p}, \]
and hence Lemma \ref{lem:Pn} gives for $|x| > 1$ and $c \in \mathbb{C}$
\begin{equation}  \label{eq:Pnasym}      
  \lim_{n\to\infty} |P_n(cx^n|q)|^{1/n^2} =  \sqrt{p}\ |x|.
\end{equation}

\begin{theorem}  \label{thm1}
Suppose $p>1$ is an integer.
Let $a_n$ and $b_n$ be given by (\ref{eq:bn})--(\ref{eq:an}), then
$a_n \in \mathbb{Z}$, $(-1)^n b_n \in \mathbb{N}$, and
$(-1)^n (b_n h_p(1) -a_n) > 0$ for $n >1$. Furthermore
\[ \lim_{n \to \infty} |b_n h_p(1) - a_n|^{1/n^2} = p^{-\frac{3(\pi^2-2)}{2\pi^2}} < 1, \]
which implies that $h_p(1)$ is irrational with measure of
irrationality $r \leq \frac{2\pi^2}{\pi^2-2} = 2.50828\ldots$
\end{theorem}
\textbf{Proof:}
If we take $x=p$ and $c=1$ in (\ref{eq:Pnasym}) then we have
\[  \lim_{n \to \infty} |P_n(p^n|q)|^{1/n^2} = p^{3/2}, \]
and combining with (\ref{eq:dn}) we have for the integer $b_n$ in (\ref{eq:bn}) 
\[  \lim_{n \to \infty}  |b_n|^{1/n^2} = p^{\frac{3(\pi^2+2)}{2\pi^2}}.  \]
Observe that $b_n$ has the same sign as $P_n(p^n|q)$ which
is $(-1)^n$. Furthermore (\ref{eq:error}) and (\ref{eq:bounds}) show that
\[  \lim_{n\to \infty} |b_n h_p(1)-a_n|^{1/n^2} = 
 \frac{\lim_{n \to \infty} d_n(p)^{1/n^2}}
 {\lim_{n\to \infty} |P_n(p^n|q)|^{1/n^2}} = p^{-\frac{3(\pi^2-2)}{2\pi^2}}<1. \]
The irrationality now follows from Lemma \ref{lem1}. Observe that
this gives rational approximants $a_n/b_n$ for $h_p(1)$ satisfying
\[   | h_p(1) - \frac{a_n}{b_n} | = 
\mathcal{O} \left(\frac{p^{(\frac{-3(\pi^2-2)}{2\pi^2}+\epsilon)n^2}}{b_n} \right) \]
for every $\epsilon > 0$. Now $b_n = p^{3n^2(\pi^2+2)/(2\pi^2) + o(n^2)}$, hence
\[   |h_p(1) - \frac{a_n}{b_n} | = \mathcal{O}
  \left( \frac{1}{b_n^{1+\frac{\pi^2-2}{\pi^2+2} - \epsilon}} \right) \]
for every $\epsilon > 0$, which gives for the measure of
irrationality the bound $r \leq 1 + \frac{\pi^2+2}{\pi^2-2} = \frac{2\pi^2}{\pi^2-2}$. \ \framebox{}  

\bigskip
The upper bound for the measure of irrationality is better than the upper bound
$4.8$ obtained in \cite{WZ}, but the same as the (earlier) upper bound of
Bundschuh and V\"a\"an\"anen \cite{BV}.
      
\section{The $q$-analog of the logarithm of 2}

Next we show that a very similar analysis also proves the irrationality of
$\ln_p(2)$ for every integer $p > 1$. First of all we rewrite
$\ln_p(2)$ using the geometric series
\[  \ln_p(2) = \sum_{k=1}^\infty (-1)^k \frac{q^k}{1-q^k} 
             = \sum_{k=1}^\infty (-1)^k q^k \sum_{j=0}^\infty q^{jk}. \]
Fubini's theorem allows us to change the order of the sums whenever
$0 < q < 1$ and this gives
\begin{eqnarray*}
  \ln_p(2) & = & \sum_{j=0}^\infty \sum_{k=1}^\infty (-1)^k q^{k(j+1)}  \\
        & = & - \sum_{j=0}^\infty \frac{q^{j+1}}{1+q^{j+1}} \\
        & = & - \sum_{k=1}^\infty \frac{1}{p^k+1}.
\end{eqnarray*}
Hence if we evaluate the Stieltjes function (\ref{eq:markov2}) at
$z=-p^n$ then we find
\[   f(-p^n) = - \sum_{k=0}^\infty \frac{1}{p^{n+k}+1}
             = \ln_p(2) + \sum_{k=1}^{n-1} \frac{1}{p^k+1}, \]
so that $f(-p^n)$ gives the required $\ln_p(2)$ up to
$\sum_{k=1}^{n-1} 1/(p^k+1)$, which is a rational number.
We can now proceed as in the previous section and evaluate
(\ref{eq:pade}) for little $q$-Legendre polynomials at $z=-p^n$
to find
\[   P_n(-p^n|q) \left( \ln_p(2) + \sum_{k=1}^{n-1} \frac{1}{p^k+1}
  \right) - Q_n(-p^n|q) = - \sum_{k=0}^\infty 
  \frac{P_n(q^k|q)}{p^n+q^k} q^k. \]             
Here we can use (\ref{eq:Pnpbinom}) to see that
\[    P_n(-p^n|q) = \sum_{k=0}^n \qbinom{n}{k}_p \qbinom{n+k}{k}_p
   p^{k(k-1)/2}, \]
is a positive integer, and if we use (\ref{eq:Qnx}) and
$(-p^nq^{\ell+1};q)_{k-\ell} =(-p^{n-k};p)_{k-\ell}$ then 
\[  Q_n(-p^n|q) = 
    (-1)^{n+1} \sum_{k=0}^n \qbinom{n}{k}_p \qbinom{n+k}{k}_p (-1)^k
                  p^{(n-k)(n-k+1)/2} \sum_{\ell=1}^{k} \frac{(-p^{n-k};p)_{k-\ell}}{p^\ell-1} . \]
This quantity has the numbers $p^\ell-1$ $(\ell=1,2,\ldots,n)$ in the denominator,
so that we need to multiply it by a multiple of all $p^\ell-1$ with $1 \leq \ell \leq n$.                                
Now we have an additional quantity
\[  \sum_{k=1}^{n-1} \frac{1}{p^k+1}  \]
and in order to make this an integer we need to multiply
it by a multiple of all $p^k+1$ for $k=1,2,\ldots,n-1$. If we choose
\[   \hat{d}_n(p) = \prod_{k=1}^n \Phi_k(p^2), \]
then because of $p^{2k}-1=(p^k-1)(p^k+1)$ we see that $\hat{d}_n(p)$
is a multiple of all $p^k+1$ and all $p^k-1$ for $k=1,\ldots,n$. Note
that $\hat{d}_n(p) = d_n(p^2)$, where $d_n$ is given by (\ref{eq:lcm}), so
that Lemma \ref{lm:dn} gives the growth
\begin{equation}   \label{eq:d2n}
   \lim_{n \to \infty} d_n(p^2)^{1/n^2} = p^{6/\pi^2}. 
\end{equation}
So if we choose
\begin{eqnarray}            
           b_n & = & d_n(p^2) P_n(-p^n|q), 
\label{eq:bnl} \\
           a_n & = & d_n(p^2) Q_n(-p^n|q) -
            b_n \sum_{k=1}^{n-1} \frac{1}{p^k+1}, \label{eq:anl}
\end{eqnarray}             
then $a_n$ and $b_n$ are integers and
\begin{equation}  \label{eq:errorl}
   b_n \ln_p(2) - a_n = d_n(p^2) \sum_{k=0}^\infty 
   \frac{P_n(q^k|q)}{-p^n-q^j} q^k .
\end{equation} 

\begin{theorem}  \label{thm2}
Suppose $p>1$ is an integer.
Let $a_n$ and $b_n$ be given by (\ref{eq:bnl})--(\ref{eq:anl}), then
$a_n \in \mathbb{Z}$, $b_n \in \mathbb{N}$, and $b_n \ln_p(2) - a_n < 0$.
Furthermore
\[  \lim_{n \to \infty} |b_n \ln_p(2) - a_n |^{1/n^2} = p^{-\frac{3(\pi^2-4)}{2\pi^2}} < 1, \] 
which implies that $\ln_p(2)$ is irrational. Its measure of irrationality
satisfies $r \leq \frac{2\pi^2}{\pi^2-4} = 3.36295\ldots$.
\end{theorem}
\textbf{Proof:}
Use $c=-1$ and $x=p$ in (\ref{eq:Pnasym}), then together with (\ref{eq:d2n}) we find
\[  \lim_{n \to \infty} b_n^{1/n^2} = p^{\frac{3(\pi^2+4)}{2\pi^2} n^2}.  \] 
Furthermore, we have
\[   \sum_{k=0}^\infty 
   \frac{P_n(q^k|q)}{-p^n-q^j} q^k = \frac{1}{P_n(-p^n|q)}
   \sum_{k=0}^\infty 
   \frac{P_n^2(q^k|q)}{-p^n-q^j} q^k \]
and $p^n \leq p^n+q^j \leq p^n+1$ for every $j$, combined with (\ref{eq:errorl}) implies
\begin{equation*}
  \frac{d_n(p^2)}{P_n(-p^n|q)}
   \frac{p^{n+1}}{(p^n+1)(p^{2n+1}-1)} \leq
  -(b_n \ln_p(2) -a_n) 
\leq \frac{d_n(p^2)}{P_n(-p^n|q)}
    \frac{p}{p^{2n+1}-1}. 
\end{equation*}
 From this one easily finds the required asymptotics and the irrationality 
then follows from Lemma \ref{lem1}. Observe that $b_n = p^{3n^2(\pi^2+4)/(2\pi^2)+o(n^2)}$
and
\[   | \ln_p(2) - \frac{a_n}{b_n} | = \mathcal{O} \left(
   \frac{1}{b_n^{1+\frac{\pi^2-4}{\pi^2+4}-\epsilon}} \right)  \]
 for every $\epsilon > 0$, which gives for the measure of
 irrationality the bound $r \leq 1+\frac{\pi^2+4}{\pi^2-4}=\frac{2\pi^2}{\pi^2-4}$. \ \framebox{}
\bigskip

The upper bound is better than the upper bounds 4.8 (obtained in \cite{WZ}), 4.311 (obtained in \cite{BV}),
 and 3.9461 (obtained in \cite{MV}).      
   
\section{Extensions}           
The construction of rational approximants for $h_p(1)$ and $\ln_p(2)$
can be extended with little effort to series of the form
\[    L = \sum_{k=1}^\infty \frac{1}{c p^k-1}, \]
where $c = a/b$ is a rational number and $c p^k \neq 1$ for every $k \geq 1$.
Indeed, these series can be obtained by evaluating the Stieltjes
function $f$ in (\ref{eq:markov2}) at $c p^n$, giving
\[  f(cp^n) = \sum_{k=1}^\infty \frac{1}{cp^k-1} -
   \sum_{k=1}^{n-1} \frac{1}{cp^k-1}. \]
We then get
\[   P_n(cp^n|q) \left( L - \sum_{k=1}^{n-1} \frac{1}{cp^k-1}
  \right) - Q_n(cp^n|q) =  \sum_{k=0}^\infty 
  \frac{P_n(q^k|q)}{cp^n-q^k} q^k, \]
where
\[  P_n(cp^n|q) =   \sum_{k=0}^n \qbinom{n}{k}_p \qbinom{n+k}{k}_p
    p^{k(k-1)/2} (-c)^k , \]
and
\[ Q_n(cp^n|q) = (-1)^{n+1} \sum_{k=0}^n \qbinom{n}{k}_p \qbinom{n+k}{k}_p (-1)^k
                  p^{(n-k)(n-k+1)/2} \sum_{\ell=1}^{k} \frac{(cp^{n-k};p)_{k-\ell}}{p^\ell-1} . \] 
 In order to have integers,
the quantities $P_n(cp^n|q)$ and $Q(cp^n|q)$ now need to be multiplied by 
$b^n$ and by a multiple of all $p^\ell-1$ for $1\leq \ell \leq n$ and of
all $ap^k-b$ for $1\leq k \leq n-1$. A possible factor is
$b^{2n} d_n(p) (c;p)_n$. This factor grows like
\[  \lim_{n \to \infty}  |b^{2n} d_n(p) (c;p)_n|^{1/n^2} = p^{\frac{3}{\pi^2} + \frac12}. \]  
In a way similar to the proof of Theorems \ref{thm1} and \ref{thm2} we can then prove:

\begin{theorem}  \label{thm3}
Suppose $p>1$ is an integer and $c = a/b$ is rational but $cp^k \neq 1$ for every $k=1,2,\ldots$
Let $a_n$ and $b_n$ be given by 
\begin{eqnarray}
    b_n & = & b^{2n} d_n(p) (c;p)_n P_n(cp^n|q), \\
    a_n & = & b^{2n} d_n(p) (c;p)_n Q_n(cp^n|q) +b_n \sum_{k=1}^{n-1} \frac{1}{cp^k-1}.
\end{eqnarray}
Then
$a_n, b_n \in \mathbb{Z}$, and
\[  \lim_{n \to \infty} |b_n L - a_n |^{1/n^2}
 = p^{-\frac{\pi^2-3}{\pi^2}} < 1, \] 
which implies that the infinite sum $L$ is irrational. Its measure of irrationality
satisfies $r \leq \frac{3\pi^2}{\pi^2-3} = 4.310119\ldots$.
\end{theorem}

The upper bound for the measure of irrationality corresponds to the upper bound
given by Bundschuh and V\"a\"an\"anen
 \cite[p.~178]{BV}. For the cases $c=1$ and $c=-1$, which we handled in
Theorems \ref{thm1} and \ref{thm2}, one can find better upper bounds.
Note that the results in \cite{BV} and \cite{MV} are also valid for
$p$ and $c$ in other number fields. Our main purpose in this paper, however,
was to emphasise the use of little $q$-Legendre polynomials in
the construction of rational approximants for certain important Lambert
series.

\section*{Acknowledgments}
This research was carried out while visiting Georgia
Institute of Technology. The author wishes to thank the School of
Mathematics for its hospitality. Also thanks to a referee for
pointing out references \cite{Be}, \cite{BV}, \cite{Ge}, and \cite{MV},
and to Roberto Costas for a useful conversation that led to Lemma
\ref{lm:dn}.  This research is partially funded by
FWO research project G.0278.97 and INTAS 2000-272.

{\obeylines
Department of Mathematics
Katholieke Universiteit Leuven
Celestijnenlaan 200 B
B-3001 Leuven
BELGIUM
walter@wis.kuleuven.ac.be}
\end{document}